\documentclass[oneside,english]{amsart}
\usepackage[T1]{fontenc}
\usepackage[latin9]{inputenc}
\usepackage{babel}
\usepackage{amstext}
\usepackage{amsthm}
\usepackage[unicode=true,pdfusetitle,
 bookmarks=true,bookmarksnumbered=false,bookmarksopen=false,
 breaklinks=false,pdfborder={0 0 1},backref=false,colorlinks=false]
 {hyperref}

\makeatletter
\numberwithin{equation}{section}
\theoremstyle{plain}
\newtheorem{thm}{\protect\theoremname}[section]

\makeatother

\providecommand{\theoremname}{Theorem}

\begin{document}
\title{Some truncated identities of Gauss}
\author{Bing He}
\address{School of Mathematics and Statistics, Central South University \\
Changsha 410083, Hunan, People's Republic of China}
\email{yuhe001@foxmail.com; yuhelingyun@foxmail.com}
\begin{abstract}
Motivated by the works of Andrews\textendash Merca and Guo\textendash Zeng,
we establish some truncated identities of Gauss by using some summation
formulas from the works of Zhi-Guo Liu \cite{L1,L6}. These give three
new expansions for partial sums of Gauss\textquoteright{} triangular
series.
\end{abstract}

\maketitle

\section{Introduction}

Euler\textquoteright s pentagonal number theorem \cite[Corollary  1.7]{A98}
may be one of the most important formulas in basic hypergeometric
series, which states that 
\[
(q;q)_{\infty}=\sum_{n=-\infty}^{\infty}(-1)^{n}q^{n(3n-1)/2},
\]
where
\[
(z;q)_{\infty}:=\prod_{n=0}^{\infty}(1-zq^{n}).
\]
Here and in the sequel, we always assume that $|q|<1$. Using induction,
Andrews and Merca \cite{AM} established the following truncated version
for Euler\textquoteright s pentagonal number theorem: 
\[
\frac{1}{(q;q)_{\infty}}\sum_{j=0}^{n-1}(-1)^{j}q^{j(3j+1)/2}(1-q^{2j+1})=1+(-1)^{n-1}\sum_{k=1}^{\infty}\frac{q^{{n \choose 2}+(n+1)k}}{(q;q)_{k}}{k-1 \brack n-1}_{q},
\]
where 
\[
(a;q)_{l}:=\frac{(a;q)_{\infty}}{(aq^{l};q)_{\infty}}
\]
and 
\[
{M \brack N}_{q}:=\begin{cases}
\frac{(q;q)_{M}}{(q;q)_{N}(q;q)_{M-N}}, & \mathit{if}\;0\leq N\leq M,\\
0, & \mathit{otherwise}.
\end{cases}
\]
Apart from the Euler pentagonal number theorem, there are two other
classical theta identities due to Gauss \cite[Corollary 2.10]{A98}:
\begin{align}
 & \frac{(q;q)_{\infty}}{(-q;q)_{\infty}}=\sum_{n=-\infty}^{\infty}(-1)^{n}q^{n^{2}},\label{eq:0-1}\\
 & \frac{(q^{2};q^{2})_{\infty}}{(q;q^{2})_{\infty}}=\sum_{n=0}^{\infty}q^{n(n+1)/2}.\label{eq:0-2}
\end{align}
With series munipulations we can show that \eqref{eq:0-2} is equivalent
to 
\begin{equation}
\frac{(q^{2};q^{2})_{\infty}}{(-q;q^{2})_{\infty}}=\sum_{n=-\infty}^{\infty}(-1)^{n}q^{n(2n+1)}.\label{eq:0-7}
\end{equation}
and 
\begin{equation}
\frac{(q^{2};q^{2})_{\infty}}{(-q;q^{2})_{\infty}}=\sum_{j=0}^{\infty}(-1)^{j}q^{j(2j+1)}(1-q^{2j+1}).\label{eq:0-6}
\end{equation}

Motivated by the work of Andrews and Merca, Guo and Zeng \cite{GZ}
showed the following truncated forms for \eqref{eq:0-1} and \eqref{eq:0-6}:
\begin{equation}
\begin{aligned} & \frac{(-q;q)_{\infty}}{(q;q)_{\infty}}\sum_{j=-n}^{n}(-1)^{j}q^{j^{2}}\\
 & =1+(-1)^{n}\sum_{k=n+1}^{\infty}\frac{(-q;q)_{n}(-1;q)_{k-n}q^{(n+1)k}}{(q;q)_{k}}{k-1 \brack n}_{q}
\end{aligned}
\label{eq:0-3}
\end{equation}
and 
\[
\begin{aligned} & \frac{(-q;q^{2})_{\infty}}{(q^{2};q^{2})_{\infty}}\sum_{j=0}^{n-1}(-1)^{j}q^{j(2j+1)}(1-q^{2j+1})\\
 & =1+(-1)^{n-1}\sum_{k=n}^{\infty}\frac{(-q;q^{2})_{n}(-q;q^{2})_{k-n}q^{2(n+1)k-n}}{(q^{2};q^{2})_{k}}{k-1 \brack n-1}_{q^{2}}.
\end{aligned}
\]

Our motivation for the present work emanates from the works of Andrews\textendash Merca
\cite{AM} and Guo\textendash Zeng \cite{GZ}. In this work, we shall
establish two truncated forms for \eqref{eq:0-2} and one truncated
version for \eqref{eq:0-7}.
\begin{thm}
\label{t1-1} For any nonnegative integer $n$ we have
\begin{align*}
 & \frac{(q;q^{2})_{\infty}}{(q^{2};q^{2})_{\infty}}\sum_{j=0}^{n}q^{j(j+1)/2}\\
 & =1+(-1)^{n}q^{n(n+1)/2}\sum_{l=n+1}^{\infty}(-q)^{l}\sum_{j=0}^{l}\frac{(q^{-n-1/2};q)_{j}(q^{1/2};q)_{l-j}}{(q;q)_{j}(q;q)_{l-j}}(-q^{n+1})^{j}{l-1 \brack n}_{q}.
\end{align*}
\end{thm}
\begin{thm}
\label{t1-2} For any nonnegative integer $n$ we have
\begin{align*}
 & \frac{(q;q^{2})_{\infty}}{(q^{2};q^{2})_{\infty}}\sum_{j=0}^{2n}q^{j(j+1)/2}\\
 & =1-(q;q^{2})_{n+1}\sum_{k=n+1}^{\infty}\frac{(q;q^{2})_{k-n-1}q^{(2n+2)k-n-1}}{(1-q^{2k})(q^{2};q^{2})_{k-1}}{k-1 \brack n}_{q^{2}}.
\end{align*}
\end{thm}
\begin{thm}
\label{t1-4} For any nonnegative integer $n$ we have
\begin{align*}
 & \frac{(-q;q^{2})_{\infty}}{(q^{2};q^{2})_{\infty}}\sum_{j=-n}^{n}(-1)^{j}q^{j(2j+1)}\\
 & =1+(-1)^{n}(-q;q^{2})_{n+1}\sum_{k=n+1}^{\infty}\frac{(-q;q^{2})_{k-n-1}}{(q^{2};q^{2})_{k}}q^{(n+1)(2k-1)}{k-1 \brack n}_{q^{2}}.
\end{align*}
\end{thm}
Let $\mathrm{pod}(n)$ count the number of partitions of $n$ without
repeated odd parts. Then the generating function for $\mathrm{pod}(n)$
is 
\[
\sum_{n=0}^{\infty}\mathrm{pod}(n)q^{n}=\frac{(-q;q^{2})_{\infty}}{(q^{2};q^{2})_{\infty}}.
\]
As a consequence of Theorem \ref{t1-4} we have the following positivity
theorem.
\begin{thm}
For any nonnegative integer $N$ we have
\[
(-1)^{n}\sum_{j=-n}^{n}(-1)^{j}\mathrm{pod}(N-j(2j+1))\geq0.
\]
\end{thm}
We collect some facts about basic hypergeometric series in the next
section. We can prove the formula \eqref{eq:0-3} in the same way
as Theorems \ref{t1-1} and \ref{t1-4} by adopting our method. In
Section \ref{sec:new} we provide a new proof of \eqref{eq:0-3}.
Section \ref{sec:1-2} is devoted to our proofs of Theorems \ref{t1-1}
and \ref{t1-2}. In the last section we will show Theorem \ref{t1-4}.

Lastly, it should be emphasized that the essential tools for our proofs
of the results are formulas from Liu's works \cite{L1,L6}.

\section{Preliminaries}

 In this section we collect several useful facts on basic hypergeometric
series.

Throughout this paper we use the following compact $q$-notation:
\[
(a_{1},a_{2},\cdots,a_{m};q)_{n}:=(a_{1};q)_{n}(a_{2};q)_{n}\cdots(a_{m};q)_{n},
\]
where $n$ is an integer or $\infty$.

The basic hypergeometric series $_{r+1}\phi_{r}$ is defined by \cite[(1.2.22)]{GR}
\[
_{r+1}\phi_{r}\left(\begin{matrix}a_{0},a_{1},\cdots,a_{r}\\
b_{1},\cdots,b_{r}
\end{matrix};q,z\right):=\sum_{n=0}^{\infty}\frac{(a_{0},a_{1},\cdots,a_{r};q)_{n}}{(q,b_{1},\cdots,b_{r};q)_{n}}z^{n}.
\]

The\textbf{ $q$}-binomial theorem \cite[(II.3)]{GR} is one of the
most interesting identities in basic hypergeometric series:

\begin{equation}
\sum_{n=0}^{\infty}\frac{(a;q)_{n}}{(q;q)_{n}}z^{n}=\frac{(az;q)_{\infty}}{(z;q)_{\infty}}\label{eq:1-2}
\end{equation}
where $|z|<1.$

One special case of \textbf{$q$}-binomial theorem is as follows:
\begin{equation}
\sum_{j=0}^{n}\frac{(q^{-n};q)_{j}}{(q;q)_{j}}z^{j}=(q^{-n}z;q)_{n},\label{eq:1-3}
\end{equation}
where $n$ is a nonnegative integer.

From \cite[Propositions 2.4 and 2.5]{L6}\footnote{The factor $(-1)^{n}$ is missing on the left-hand side of \cite[Propositions 2.4 and 2.5]{L6}.}
we have

\begin{equation}
\begin{aligned} & (-1)^{n}\frac{(\alpha q;q)_{n}}{(q;q)_{n}}q^{{n+1 \choose 2}}\text{}_{3}\phi_{2}\left(\begin{matrix}q^{-n},\alpha q^{n+1},\alpha cd/q\\
\alpha c,\alpha d
\end{matrix};q,1\right)\\
 & =\sum_{j=0}^{n}(-1)^{j}\frac{(1-\alpha q^{2j})(\alpha,q/c,q/d;q)_{j}}{(1-\alpha)(q,\alpha c,\alpha d;q)_{j}}q^{j(j-3)/2}(\alpha cd)^{j}
\end{aligned}
\label{eq:1-2-1}
\end{equation}
and
\begin{equation}
\begin{aligned} & (-1)^{n}\frac{(\alpha q;q)_{n}}{(q;q)_{n}}q^{{n+1 \choose 2}}\text{}_{2}\phi_{1}\left(\begin{matrix}q^{-n},\alpha q^{n+1}\\
\alpha c
\end{matrix};q,1\right)\\
 & =\sum_{j=0}^{n}\frac{(1-\alpha q^{2j})(\alpha,q/c;q)_{j}}{(1-\alpha)(q,\alpha c;q)_{j}}q^{j^{2}-j}(\alpha c)^{j}.
\end{aligned}
\label{eq:7-12}
\end{equation}

\section{\label{sec:new} A new proof of \eqref{eq:0-3}}

Recall the following identity \cite[p. 2088]{L6}:

\begin{equation}
q^{n(n+1)/2}{}_{2}\phi_{1}\left(\begin{matrix}q^{-n},q^{n+1}\\
-q
\end{matrix};q,1\right)=\sum_{j=-n}^{n}(-1)^{n+j}q^{j^{2}}.\label{eq:1-1}
\end{equation}

It follows from \eqref{eq:1-2} and \eqref{eq:1-3} that
\begin{align*}
\frac{(-q,q)_{\infty}}{(q,q)_{\infty}}\sum_{j=0}^{n}\frac{(q^{-n},q^{n+1};q)_{j}}{(q,-q;q)_{j}} & =\frac{(-q,q)_{\infty}}{(q;q)_{n}(q,q)_{\infty}}\sum_{j=0}^{n}\frac{(q^{-n};q)_{j}(q;q)_{n+j}}{(q,-q;q)_{j}}\\
 & =\frac{1}{(q;q)_{n}}\sum_{j=0}^{n}\frac{(q^{-n};q)_{j}}{(q;q)_{j}}\frac{(-q^{j+1};q)_{\infty}}{(q^{n+j+1};q)_{\infty}}\\
 & =\frac{1}{(q;q)_{n}}\sum_{j=0}^{n}\frac{(q^{-n};q)_{j}}{(q;q)_{j}}\sum_{k=0}^{\infty}\frac{(-q^{-n};q)_{k}}{(q;q)_{k}}q^{(n+j+1)k}\\
 & =\frac{1}{(q;q)_{n}}\sum_{k=0}^{\infty}\frac{(-q^{-n};q)_{k}}{(q;q)_{k}}q^{(n+1)k}\sum_{j=0}^{n}\frac{(q^{-n};q)_{j}}{(q;q)_{j}}q^{kj}\\
 & =\frac{1}{(q;q)_{n}}\sum_{k=0}^{\infty}\frac{(-q^{-n};q)_{k}(q^{k-n};q)_{n}}{(q;q)_{k}}q^{(n+1)k}.
\end{align*}
Then, by \eqref{eq:1-1},
\begin{align*}
 & \frac{(-q,q)_{\infty}}{(q,q)_{\infty}}\sum_{j=-n}^{n}(-1)^{j}q^{j^{2}}\\
 & =(-1)^{n}q^{n(n+1)/2}\frac{(-q,q)_{\infty}}{(q,q)_{\infty}}\sum_{j=0}^{n}\frac{(q^{-n},q^{n+1};q)_{j}}{(q,-q;q)_{j}}\\
 & =\frac{(-1)^{n}q^{n(n+1)/2}}{(q;q)_{n}}\sum_{k=0}^{\infty}\frac{(-q^{-n};q)_{k}(q^{k-n};q)_{n}}{(q;q)_{k}}q^{(n+1)k}\\
 & =\frac{(-1)^{n}q^{n(n+1)/2}}{(q;q)_{n}}(q^{-n};q)_{n}+\frac{(-1)^{n}q^{n(n+1)/2}}{(q;q)_{n}}\sum_{k=n+1}^{\infty}\frac{(-q^{-n};q)_{k}(q^{k-n};q)_{n}}{(q;q)_{k}}q^{(n+1)k}\\
 & =1+(-1)^{n}(-q;q)_{n}\sum_{k=n+1}^{\infty}\frac{(-1;q)_{k-n}}{(q;q)_{k}}q^{(n+1)k}{k-1 \brack n}_{q},
\end{align*}
where for the second to the last equality we applied $(q^{k-n};q)_{n}=0$
for $0<k\leq n$ and the last equality follows from the identity $(-q^{-n};q)_{k}=q^{-n(n+1)/2}(-q;q)_{n}(-1;q)_{k-n}.$
\qed

\section{\label{sec:1-2} Proofs of Theorems \ref{t1-1} and \ref{t1-2}}

\noindent{\it Proof of Theorem \ref{t1-1}.} It follows from \eqref{eq:1-2}
that 
\begin{align*}
\frac{(q^{2k+3};q^{2})_{\infty}}{(-q^{k+1},q^{n+k+2};q)_{\infty}} & =\frac{(-q^{k+3/2},q^{k+3/2};q)_{\infty}}{(-q^{k+1},q^{n+k+2};q)_{\infty}}\\
 & =\sum_{i=0}^{\infty}\frac{(q^{1/2};q)_{i}}{(q;q)_{i}}(-q^{k+1})^{i}\sum_{j=0}^{\infty}\frac{(q^{-n-1/2};q)_{j}}{(q;q)_{j}}(q^{n+k+2})^{j}\\
 & =\sum_{l=0}^{\infty}(-q^{k+1})^{l}\sum_{j=0}^{\infty}\frac{(-1)^{j}(q^{-n-1/2};q)_{j}(q^{1/2};q)_{l-j}}{(q;q)_{j}(q;q)_{l-j}}(q^{n+1})^{j},
\end{align*}
where the last equality follows by making the change of variables
$i=l-j.$ Then, by \eqref{eq:1-3},
\begin{equation}
\begin{aligned} & \frac{(q^{3};q^{2})_{\infty}}{(q^{2};q^{2})_{\infty}}\sum_{k=0}^{n}\frac{(q^{-n},-q;q)_{k}(q^{n+1};q)_{k+1}}{(q;q)_{k}(q^{3};q^{2})_{k}}\\
 & =\frac{(q^{3};q^{2})_{\infty}}{(q;q)_{n}(q^{2};q^{2})_{\infty}}\sum_{k=0}^{n}\frac{(q^{-n},-q;q)_{k}(q;q)_{n+k+1}}{(q;q)_{k}(q^{3};q^{2})_{k}}\\
 & =\frac{1}{(q;q)_{n}}\sum_{k=0}^{n}\frac{(q^{-n};q)_{k}}{(q;q)_{k}}\frac{(q^{2k+3};q^{2})_{\infty}}{(-q^{k+1},q^{n+k+2};q)_{\infty}}\\
 & =\frac{1}{(q;q)_{n}}\sum_{l=0}^{\infty}(-q)^{l}\sum_{j=0}^{\infty}\frac{(q^{-n-1/2};q)_{j}(q^{1/2};q)_{l-j}}{(q;q)_{j}(q;q)_{l-j}}(-q^{n+1})^{j}\sum_{k=0}^{n}\frac{(q^{-n};q)_{k}}{(q;q)_{k}}q^{lk}\\
 & =\frac{1}{(q;q)_{n}}\sum_{l=0}^{\infty}(-q)^{l}\sum_{j=0}^{l}\frac{(q^{-n-1/2};q)_{j}(q^{1/2};q)_{l-j}}{(q;q)_{j}(q;q)_{l-j}}(-q^{n+1})^{j}(q^{l-n};q)_{n}.
\end{aligned}
\label{eq:4-7}
\end{equation}
Set $\alpha=q,c=q^{1/2},d=-q^{1/2}$ in \eqref{eq:1-2-1}. We obtain
\[
_{3}\phi_{2}\left(\begin{matrix}q^{-n},q^{n+2},-q\\
q^{3/2},-q^{3/2}
\end{matrix};q,1\right)=(-1)^{n}\frac{1-q}{1-q^{n+1}}q^{-n(n+1)/2}\sum_{j=0}^{n}q^{j(j+1)/2}.
\]
Then, by the last equality of \eqref{eq:4-7},
\begin{align*}
 & \frac{(q;q^{2})_{\infty}}{(q^{2};q^{2})_{\infty}}\sum_{j=0}^{n}q^{j(j+1)/2}\\
 & =(-1)^{n}q^{n(n+1)/2}\frac{(q^{3};q^{2})_{\infty}}{(q^{2};q^{2})_{\infty}}\sum_{j=0}^{n}\frac{(q^{-n},-q;q)_{j}(q^{n+1};q)_{j+1}}{(q;q)_{j}(q^{3};q^{2})_{j}}\\
 & =\frac{(-1)^{n}q^{n(n+1)/2}}{(q;q)_{n}}\sum_{l=0}^{\infty}(-q)^{l}\sum_{j=0}^{l}\frac{(q^{-n-1/2};q)_{j}(q^{1/2};q)_{l-j}}{(q;q)_{j}(q;q)_{l-j}}(-q^{n+1})^{j}(q^{l-n};q)_{n}
\end{align*}
\begin{align*}
 & =\frac{(-1)^{n}q^{n(n+1)/2}}{(q;q)_{n}}(q^{-n};q)_{n}\\
 & +\frac{(-1)^{n}q^{n(n+1)/2}}{(q;q)_{n}}\sum_{l=n+1}^{\infty}(-q)^{l}\sum_{j=0}^{l}\frac{(q^{-n-1/2};q)_{j}(q^{1/2};q)_{l-j}}{(q;q)_{j}(q;q)_{l-j}}(-q^{n+1})^{j}(q^{l-n};q)_{n}\\
 & =1+(-1)^{n}q^{n(n+1)/2}\sum_{l=n+1}^{\infty}(-q)^{l}\sum_{j=0}^{l}\frac{(q^{-n-1/2};q)_{j}(q^{1/2};q)_{l-j}}{(q;q)_{j}(q;q)_{l-j}}(-q^{n+1})^{j}{l-1 \brack n}_{q}.
\end{align*}
This completes the proof. \qed

\noindent{\it Proof of Theorem \ref{t1-2}.} Recall the following
identity \cite[(7.15)]{L1}:
\[
\sum_{j=0}^{n}\frac{(q;q^{2})_{j}q^{-n(2j+1)}}{(q^{2};q^{2})_{j}}=\frac{(q;q^{2})_{n}}{(q^{2};q^{2})_{n}}\sum_{j=0}^{2n}q^{-{j+1 \choose 2}}.
\]
Replacing $q$ by $q^{-1}$ in this identity we have
\[
\sum_{j=0}^{n}\frac{(q;q^{2})_{j}q^{j(2n+1)}}{(q^{2};q^{2})_{j}}=\frac{(q;q^{2})_{n}}{(q^{2};q^{2})_{n}}\sum_{j=0}^{2n}q^{{j+1 \choose 2}}.
\]
This formula can also be obtained by using induction on $n.$ Then
\begin{align*}
\frac{(q;q^{2})_{\infty}}{(q^{2};q^{2})_{\infty}}\sum_{j=0}^{2n}q^{j(j+1)/2} & =\frac{(q^{2n+1};q^{2})_{\infty}}{(q^{2n+2};q^{2})_{\infty}}\sum_{j=0}^{n}\frac{(q;q^{2})_{j}q^{j(2n+1)}}{(q^{2};q^{2})_{j}}\\
 & =\sum_{j=0}^{n}\frac{(q;q^{2})_{j}q^{(2n+1)j}}{(q^{2};q^{2})_{j}}\sum_{i=0}^{\infty}\frac{(q^{-1};q^{2})_{i}}{(q^{2};q^{2})_{i}}q^{(2n+2)i}\\
 & =\sum_{k=0}^{\infty}\sum_{j=0}^{n}\frac{(q;q^{2})_{j}(q^{-1};q^{2})_{k-j}q^{(2n+2)k-j}}{(q^{2};q^{2})_{j}(q^{2};q^{2})_{k-j}},
\end{align*}
where the last equality follows by making the change of variables
$i=k-j.$ Using induction on $n$ we deduce that 
\[
\sum_{j=0}^{n}\frac{(q;q^{2})_{j}(q^{-1};q^{2})_{k-j}q^{-j}}{(q^{2};q^{2})_{j}(q^{2};q^{2})_{k-j}}=-\frac{(q;q^{2})_{n+1}(q;q^{2})_{k-n-1}q^{-n-1}}{(1-q^{2k})(q^{2};q^{2})_{k-n-1}(q^{2};q^{2})_{n}}
\]
for $k\geq1.$ Thus,
\begin{align*}
 & \frac{(q;q^{2})_{\infty}}{(q^{2};q^{2})_{\infty}}\sum_{j=0}^{2n}q^{j(j+1)/2}\\
 & =1-\frac{(q;q^{2})_{n+1}}{(q^{2};q^{2})_{n}}\sum_{k=1}^{\infty}\frac{(q;q^{2})_{k-n-1}q^{(2n+2)k-n-1}}{(1-q^{2k})(q^{2};q^{2})_{k-n-1}}\\
 & =1-(q;q^{2})_{n+1}\sum_{k=n+1}^{\infty}\frac{(q;q^{2})_{k-n-1}q^{(2n+2)k-n-1}}{(1-q^{2k})(q^{2};q^{2})_{k-1}}{k-1 \brack n}_{q^{2}}.
\end{align*}
This finishes the proof. \qed

\section{Proof of Theorem \ref{t1-4}}

By \eqref{eq:1-2} and \eqref{eq:1-3} we have
\begin{equation}
\begin{aligned} & \frac{(-q;q^{2})_{\infty}}{(q^{2};q^{2})_{\infty}}\sum_{j=0}^{n}\frac{(q^{-2n},q^{2n+2};q^{2})_{j}}{(q^{2},-q;q^{2})_{j}}\\
 & =\frac{(-q;q^{2})_{\infty}}{(q^{2};q^{2})_{n}(q^{2};q^{2})_{\infty}}\sum_{j=0}^{n}\frac{(q^{-2n};q^{2})_{j}(q^{2};q^{2})_{n+j}}{(q^{2},-q;q^{2})_{j}}\\
 & =\frac{1}{(q^{2};q^{2})_{n}}\sum_{j=0}^{n}\frac{(q^{-2n};q^{2})_{j}}{(q^{2};q^{2})_{j}}\frac{(-q^{2j+1};q^{2})_{\infty}}{(q^{2(n+j+1)};q^{2})_{\infty}}\\
 & =\frac{1}{(q^{2};q^{2})_{n}}\sum_{j=0}^{n}\frac{(q^{-2n};q^{2})_{j}}{(q^{2};q^{2})_{j}}\sum_{k=0}^{\infty}\frac{(-q^{-2n-1};q^{2})_{k}}{(q^{2};q^{2})_{k}}q^{2(n+j+1)k}\\
 & =\frac{1}{(q^{2};q^{2})_{n}}\sum_{k=0}^{\infty}\frac{(-q^{-2n-1};q^{2})_{k}}{(q^{2};q^{2})_{k}}q^{2(n+1)k}\sum_{j=0}^{n}\frac{(q^{-2n};q^{2})_{j}}{(q^{2};q^{2})_{j}}q^{2kj}\\
 & =\frac{1}{(q^{2};q^{2})_{n}}\sum_{k=0}^{\infty}\frac{(-q^{-2n-1};q^{2})_{k}(q^{2k-2n};q^{2})_{n}}{(q^{2};q^{2})_{k}}q^{2(n+1)k}.
\end{aligned}
\label{eq:5-1}
\end{equation}

Replacing $q$ by $q^{2}$ and then setting $\alpha=1,c=-q$ in \eqref{eq:7-12}
we get
\[
(-1)^{n}q^{n(n+1)}{}_{2}\phi_{1}\left(\begin{matrix}q^{-2n},q^{2n+2}\\
-q
\end{matrix};q^{2},1\right)=\sum_{j=-n}^{n}(-1)^{j}q^{2j^{2}+j}.
\]
Then, by the last equality of \eqref{eq:5-1}, 
\begin{align*}
 & \frac{(-q;q^{2})_{\infty}}{(q^{2};q^{2})_{\infty}}\sum_{j=-n}^{n}(-1)^{j}q^{2j^{2}+j}\\
 & =(-1)^{n}q^{n(n+1)}\frac{(-q;q^{2})_{\infty}}{(q^{2};q^{2})_{\infty}}\sum_{j=0}^{n}\frac{(q^{-2n},q^{2n+2};q^{2})_{j}}{(q^{2},-q;q^{2})_{j}}\\
 & =\frac{(-1)^{n}q^{n(n+1)}}{(q^{2};q^{2})_{n}}\sum_{k=0}^{\infty}\frac{(-q^{-2n-1};q^{2})_{k}(q^{2k-2n};q^{2})_{n}}{(q^{2};q^{2})_{k}}q^{2(n+1)k}\\
 & =\frac{(-1)^{n}q^{n(n+1)}}{(q^{2};q^{2})_{n}}(q^{-2n};q^{2})_{n}\\
 & \;+\frac{(-1)^{n}q^{n(n+1)}}{(q^{2};q^{2})_{n}}\sum_{k=n+1}^{\infty}\frac{(-q^{-2n-1};q^{2})_{k}(q^{2k-2n};q^{2})_{n}}{(q^{2};q^{2})_{k}}q^{2(n+1)k}\\
 & =1+(-1)^{n}(-q;q^{2})_{n+1}\sum_{k=n+1}^{\infty}\frac{(-q;q^{2})_{k-n-1}}{(q^{2};q^{2})_{k}}q^{(n+1)(2k-1)}{k-1 \brack n}_{q^{2}}.
\end{align*}
This concludes the proof.

\section*{Acknowledgement}

 This work was partially supported by the National Natural Science
Foundation of China (Grant No. 11801451) and the Natural Science Foundation
of Hunan Province (Grant No. 2020JJ5682).

\end{document}